\documentclass{article}
\usepackage{amssymb,amsfonts}

\newcommand{\qed}{\hfill$\square$\par\vskip12pt}
\newcommand{\twomat}[4]{{\left(\begin{array}{cc}#1&#2\\#3&#4\end{array}\right)}}
\newcommand{\trace}{\mathop{\rm Tr}\nolimits}

\newcommand{\schatten}[2]{\left|\left|\,{#2}\,\right|\right|_{#1}}
\newcommand{\diag}{\mathop{\rm diag}\nolimits}
\newcommand{\cirk}{\mathop{\rm circ}\nolimits}

\newcommand{\cB}{{\cal B}}
\newcommand{\cC}{{\cal C}}

\newcommand{\cH}{{\cal H}}

\newcommand{\cP}{{\cal P}}
\newcommand{\cS}{{\cal S}}

\newcommand{\C}{{\mathbb{C}}}
\newcommand{\R}{{\mathbb{R}}}
\newcommand{\N}{{\mathbb{N}}}
\newcommand{\M}{{\mathbb{M}}}
\renewcommand{\u}{{\uparrow}}
\renewcommand{\d}{{\downarrow}}

\newcommand{\identity}{\mathbb{I}}
\newcommand{\id}{\identity}

\newcommand{\be}{\begin{equation}}
\newcommand{\ee}{\end{equation}}
\newcommand{\bea}{\begin{eqnarray}}
\newcommand{\eea}{\end{eqnarray}}
\newcommand{\beas}{\begin{eqnarray*}}
\newcommand{\eeas}{\end{eqnarray*}}
\newcommand{\text}{\mbox}
\newtheorem{theorem}{Theorem}
\newtheorem{conjecture}{Conjecture}
\newtheorem{problem}{Problem}
\newtheorem{proposition}{Proposition}

\begin{document}


\title{Problems and Conjectures in Matrix and Operator Inequalities}
\author{Koenraad M.R.\ Audenaert\\
Department of Mathematics, {Royal Holloway, University of London},\\
Egham TW20 0EX, United Kingdom\\[5mm]
Fuad Kittaneh\\
Department of Mathematics, University of Jordan,\\
Amman, Jordan}
\date{\today}

\maketitle

\begin{abstract}
We present a treasure trove of open problems in matrix and operator inequalities,
of a functional analytic nature, and
with various degrees of hardness.
\end{abstract}

\section{Introduction}
In this chapter we present a number of conjectures and open problems in the theory of matrix and
operator inequalities.

We have not endeavoured completeness, as doing so would require an encyclopedic volume.
Instead we have made a selection of problems that matches our research interests;
indeed, many of these problems came up in our own research work. The underlying theme
in all problems presented here is functional analysis,
with the fractional power function $x^q$ (whether or not under the guise of the Schatten norm $||\cdot||_q$) being
most prominent.
Many problems would have immediate applications in other fields of science, and we briefly indicate that.

Quite likely, the level of difficulty of these problems varies greatly. However,
while we might have included a rough indication of hardness with every problem, we thought it wiser not to.
Experience has it that the actual hardness of a problem depends on whether one starts on the right foot; that is,
with the right background, the right set of tools, the right way of looking at it.
Grading the problems according to hardness has the inherent danger of discouraging readers
and may ultimately say more about the grader than about the problems themselves.

We have randomly subdivided our list of problems into several sections, and the reader is invited to read them in any
order preferred. We hope that our readership finds these problems intellectually stimulating.

We thank Chi-Kwong Li for encouragement and suggestions, and
especially for his contribution to the material on the geometry of polynomials.
We also thank J.C.\ Bourin for his many comments.
\section{Matrix subadditivity inequalities}

If\ $f:[0,\infty )\rightarrow \lbrack 0,\infty )$ is a concave function,
then the subadditivity relation $f(a+b)\leq f(a)+f(b)$ holds for all $a,b\geq 0$. A non-commutative version
of this inequality is true for all positive semidefinite matrices.
\begin{theorem}
Let $f:[0,\infty )\rightarrow \lbrack 0,\infty )$ be concave.
If $A$ and $B$ are positive semidefinite matrices, then
for any unitarily invariant norm $||| \cdot |||$
\be
||| f(A+B)||| \leq ||| f(A)+f(B)|||.\label{eq:F}
\ee
\end{theorem}
This result is due to
Bourin and Uchiyama \cite{BU}, and for the spectral norm $|| \cdot ||$, it is due to Kosem \cite{K}.

Bourin showed \cite{B10} that inequality (\ref{eq:F}) can
be generalized to normal matrices.
\begin{theorem}[Bourin]
\label{th:B}
Let $A$ and $B$ be normal matrices and let $f:[0,\infty )\rightarrow \lbrack 0,\infty )$ be concave. Then%
\be
||| f(| A+B| )||| \leq ||| f(|A|)+f(|B| )|||
\ee
for every unitarily invariant norm.
\end{theorem}

If $A$ and/or $B$ are non-normal this inequality no longer holds, but one can ask the following:
\begin{problem}
For a given unitarily invariant norm is there is a constant $c$ such that
\be
||| f(| A+B| )||| \leq c\; ||| f(|A|)+f(|B| )||| \label{eq:Fc}
\ee
holds for all $A,B\in\M_n(\C)$ and for all concave functions $f:[0,\infty )\rightarrow \lbrack 0,\infty )$.
\end{problem}
Lee asked this question in \cite{lee11} for the special case of Schatten $p$-norms and conjectured that inequality
(\ref{eq:Fc}) holds with constant
$$
c=\sqrt{\frac{1+\sqrt{2}}{2}}
$$
for all non-negative concave $f$ when the norm is the Frobenius norm ($p=2$).
She was able to prove that (\ref{eq:Fc}) holds with $c=\sqrt{2}$ for \textit{all} unitarily invariant norms.
In that setting this constant is the best one possible as can be seen from taking the operator norm
and $f(x)=x$.

Specialising Theorem \ref{th:B} to fractional powers, we have the inequality
\be
||| ~| A+B| ^{p}||| \leq ||| ~|A| ^{p}+|B| ^{p}||| \text{ for }0<p\leq 1\text{.}
\ee
This observation prompted a number of related questions \cite{B}:
\begin{problem}[Bourin]
Given $A,B\geq 0$ and $p,q>0$, is it true that%
\be
||| ~A^{p+q}+B^{p+q}||| \leq |||(A^{p}+B^{p})(A^{q}+B^{q})||| ?
\ee
\end{problem}

\begin{problem}[Bourin]
Given $A,B\geq 0$ and $p,q>0$, is it true that%
\be
||| ~A^{p}B^{q}+B^{p}A^{q}||| \leq |||A^{p+q}+B^{p+q}||| ?
\ee
\end{problem}
W.l.o.g.\ one can take $0\le p\le 1$ and $q=1-p$ in both problems,
by absorbing $(p+q)$-th powers into $A$ and $B$.

The related inequality for Heinz means
\be
||| ~A^{p}B^{1-p}+A^{1-p}B^{p}||| \leq |||A+B|||
\ee
is known to be true.
In fact a number of stronger inequalities holds.
On one hand, for $X$ any matrix and $0\le p\le 1$ and any unitarily invariant norm,
\be
||| A^p X B^{1-p} + A^{1-p} X B^p ||| \le ||| AX+XB |||;
\ee
see e.g.\ Corollary IX.4.10 in \cite{bhatia}.
On the other hand, it was proven in \cite{ka06} that the inequality holds for each singular value:
\be
\sigma_i(A^{p}B^{1-p}+A^{1-p}B^{p}) \le \sigma_i(A+B).
\ee
Hence, the inequality of Problem 2 may similarly have a stronger counterpart for singular values:

\begin{problem}
Given $A,B\geq 0$ and $p>0$, is it true that
\be
\sigma_i(A^{p}B^{1-p}+B^{p}A^{1-p}) \le \sigma_i(A+B)?
\ee
\end{problem}
The inequality of Problem 1 has no immediate version for singular values. Counterexamples to the
inequality
$\sigma_i(A+B)\leq\sigma_i((A^{p}+B^{p})(A^{1-p}+B^{1-p}))$
are easily found.


Inequality (\ref{eq:F}) can be generalised in further ways.
Aujla and Bourin showed \cite{AB} that for monotone concave functions $f$ on $[0,+\infty)$ with $f(0)\ge0$,
for all positive semidefinite $A$ and $B$,
there exist unitary matrices $U$ and $V$ such that
\be
f(A+B)\le Uf(A)U^* +Vf(B)V^*.\label{eq:FUV}
\ee
This leads to the following questions, which were also considered in \cite{bhl}:
\begin{problem}
Can the condition of monotonicity be dropped?
That is, does (\ref{eq:FUV}) hold for all concave functions $f$ on $[0,+\infty)$ with $f(0)\ge0$?
\end{problem}
\begin{problem}
Can (\ref{eq:FUV}) be strengthened to a double inequality?
That is, under the various mentioned conditions, do there always exist unitary matrices $U$ and $V$ such that
\be
0\le f(A+B)- Uf(A)U^* \le Vf(B)V^*?\label{eq:FUV2}
\ee
\end{problem}
\section{Complementary McCarthy inequalities\label{sec:compl}}
For any known inequality between two quantities $A$ and $B$,
 $A\le B$,
it is possible to ask for complementary inequalities, to bound, say, the difference
$B-A$ from above in terms of a third quantity $C$.

A classical way to do so yields the so-called Kantorovich-type inequalities,
where $C$ involves
lower and upper bounds on the matrices appearing in $A$ and $B$.
Kantorovich's inequality itself (see e.g.
\cite{HJI}, Th. 7.4.41) states
that for any positive definite matrix $A$ with smallest eigenvalue $m$ and largest
eigenvalue $M$, and for any normalised vector $x$, we have
\be
(x,Ax)(x,A^{-1}x) \le \frac{(m+M)^2}{4mM},
\ee
which can be seen as being complementary to the inequality
\be
1\le (x,Ax)(x,A^{-1}x).
\ee
One can, however, consider complementary inequalities of a different kind, where the bounding
quantity $C$ is related
to the bounded ones, $A$ and $B$, in a more geometrical sense.
See, for example \cite{ka2}, where an inequality complementary to the
Araki-Lieb-Thirring inequality was proven.

In this section we consider complementary McCarthy inequalities.
The classic McCarthy inequality states that, for positive operators $A$ and $B$,
\be
\trace (A+B)^q \ge \trace A^q + \trace B^q,
\ee
for $q\ge1$, whereas the reversed inequality holds for $0\le q\le 1$.

Here we consider a number of upper bounds on the difference between both sides, in terms of the
third ingredient $\trace C^q$, where $C$ is related to $A$ and $B$ in a geometric way.
Clarkson's inequalities (see below) are well-known examples, where $C$ is taken to be $|A-B|$.

Another natural candidate is the geometric mean
$$
C=A\#B = A^{1/2}(A^{-1/2}BA^{-1/2})^{1/2}A^{1/2}.
$$
However, numerical simulations
indicate that the value of $\trace (A\#B)^{q}$ can be arbitrarily small for finite
values of the quantity $\trace(A+B)^q - \trace A^q-\trace B^q$.

In \cite{katrace} the related case $C=(A^{1/2}BA^{1/2})^{1/2}$ has been studied, and the following has been shown:
\begin{theorem}
For positive $A$ and $B$,
\be
\trace(A+B)^q \le \trace A^q+\trace B^q+(2^q-2)\trace (A^{1/2}BA^{1/2})^{q/2},
\ee
for $q\le -2$ and for $1\le q\le2$.
For $0<q\le 1$ and $2\le q\le3$, the reversed inequality holds.
\end{theorem}
For $1\le q\le2$ this complements McCarthy's inequality; for $2\le q\le 3$ it does not complement
but actually strengthens it.

The proof relies on several ingredients.
The first is an integral representation of negative fractional powers $x^q$ in terms of the exponential function,
\be
x^q = \frac{1}{\Gamma(-q)}\int_0^\infty \exp(-xt)\,t^{-q-1}\,dt, \quad q<0.
\ee
Simple integrations over $x$ lead to similar integrals for $x^q$ in the ranges $k<q<k+1$ $(k\in\N)$
in terms of the exponential function and polynomials of degree $k$:
\be
x^q = \frac{q}{\Gamma(1-q)} \int_0^\infty (-1)^{k+1}\left(\exp(-xt)-\sum_{j=0}^k \frac{(-xt)^j}{j!}\right)\,
t^{-q-1}\,dt,\quad k<q<k+1.
\ee
The theorem is then a relatively easy consequence of the Araki-Lieb-Thirring inequality,
the Golden-Thompson inequality for the exponential function
and the identity
\be
g(a+b)-g(a)-g(b)=g(2\sqrt{ab})-2g(\sqrt{ab}),\label{eq:gq}
\ee
which holds when $g$ is a quadratic polynomial.

The method of proof fails for $q>3$ as the identity (\ref{eq:gq}) no longer holds for polynomials of degree higher than 2.
It also fails for $-2<q<0$ as for this range the direction in which the Araki-Lieb-Thirring inequality holds is not the required one.
Nevertheless, numerical simulations heavily support the following:
\begin{conjecture}
For positive $A$ and $B$,
\be
\trace(A+B)^q \le \trace A^q+\trace B^q+(2^q-2)\trace (A^{1/2}BA^{1/2})^{q/2},
\ee
for $-2<q<0$.
For $3\le q$ the reversed inequality holds.
\end{conjecture}

\section{Matrix means}
For real scalars $x$ and $y$, many means have been defined; the simplest are the
arithmetic mean $(x+y)/2$, the geometric mean
$\sqrt{xy}$ (for non-negative $x$ and $y$) and the harmonic mean $1/(1/x+1/y)$.
We will denote a general mean by
$\mu(x,y)$. All scalar means share the property
that for all $x\le y$, $x\le \mu(x,y)\le y$. In \cite{ka_ilas} this
property is called the \textit{in-betweenness} property.

Most scalar means have been generalised to matrices and operators.
A well-known class of operator means are the Kubo-Ando means \cite{kuboando}, which include
among other means the arithmetic, geometric and harmonic mean. The Kubo-Ando means are defined
via a small number of axioms. Other matrix means exist, however, that do not belong to this
class. For example, the power means $((X^p+Y^p)/2)^{1/p}$ introduced by Bhagwat and
Subramanian \cite{bhagwat}.

Because the space of all matrices cannot be totally ordered
it is no longer clear whether or how the in-betweenness property of scalar means
can be carried over to matrix means. In \cite{ka_ilas} a first attempt was made at doing so.
The idea was to endow the matrix space with a certain metric $\delta$, and define
in-betweenness of a matrix mean $\mu(X,Y)$ as the property that
for all $X$ and $Y$, $\delta(\mu(X,Y),Y) \le \delta(X,Y)$.

Most matrix means can be generalised to weighted means by introduction of a weighting
parameter $t$; for example, the weighted power means are defined as
$(t X^p+(1-t)Y^p)^{1/p}$. We denote a weighted matrix mean by the general symbol
$\mu_t(X,Y)$.
More generally, then, than the in-betweenness property, one can ask
whether a weighted matrix mean satisfies the following monotonicity
property w.r.t. a given metric: for any $X$ and $Y$, the function
$t\mapsto \delta(\mu_t(X,Y),Y)$ is monotonically increasing over the interval
$[0,1]$.

It was shown in \cite{ka_ilas} that the weighted Kubo-Ando means are monotonic w.r.t.\
the trace metric $\delta(X,Y) = ||\log(A^{-1/2} B A^{-1/2})||_2$, but not w.r.t.\
the Euclidean metric $\delta(X,Y)=||X-Y||_2$. In contrast,
the weighted power means are monotonic w.r.t.\ the Euclidean metric for $1\le p\le2$.
There is numerical evidence that this holds for larger values of $p$ as well:
\begin{conjecture}
Let $X$ and $Y$ be positive operators, and $0\le t\le 1$.
Then the function $t\mapsto \trace(Y-(t X^p+(1-t)Y^p)^{1/p})^2$ increases
monotonically for $p\ge 2$.
\end{conjecture}
Many further questions can be raised. For example,
are the power means monotonic in other metrics as well, like the $L_q$ metric
$\delta(X,Y) = ||X-Y||_q$?
\section{Norm inequalities for partitioned operators}
Given an operator $T$ acting on a space $\cH$ that is the direct sum of $n$ spaces
$\cH_i$, $\cH=\bigoplus_i \cH_i$,
it is possible to consider the partitioned operator $[T_{ij}]$, where $T_{ij}$
is a mapping from $\cH_i$ to $\cH_j$.
In finite dimensions, a partitioned operator can be represented by a block matrix.

Given such a partitioning one can ask a variety of questions regarding the
relations between the norm of $T$ and the norms of its constituents $T_{ij}$.
One of the first papers to address this question is \cite{bk90}, where the Schatten
$q$-norm was considered. The following inequality was proven:
\be
n^{2-q} ||T||_q^q \ge \sum_{i,j=1}^n ||T_{ij}||_q^q \ge ||T||_q^q,
\ee
for $1\le q\le 2$, while for $q\ge2$ the reversed inequalities hold.

Inequalities of this kind have been called norm compression inequalities, because the
full information contained in the operator is compressed into a smaller set of quantities
(the norms of the blocks)
and the inequalities give useful bounds on the norm
of the full operator when only its compression is known.
An overview of several such inequalities can
be found in \cite{ka05} and references therein.

In \cite{ka07}  a particular kind of norm compression inequality was considered
leading to the following conjecture:
\begin{conjecture}\label{conj_ka1}
Let $T$ be the partitioned operator
$$
T = \left(
\begin{array}{cccc}
A_1 & A_2 & \cdots & A_N \\
B_1 & B_2 & \cdots & B_N
\end{array}
\right).
$$
Denoting by $\cC_p(T)$ its Schatten $p$-norm compression,
$$
\cC_p(T) =
\left(
\begin{array}{cccc}
||A_1||_p & ||A_2||_p & \cdots & ||A_N||_p \\
||B_1||_p & ||B_2||_p & \cdots & ||B_N||_p
\end{array}
\right),
$$
we have
\bea
||T||_p &\ge& \schatten{p}{\cC_p(T)},\quad 1\le p\le2 \label{eq:thedual}\\
||T||_p &\le& \schatten{p}{\cC_p(T)},\quad p\ge2. \label{eq:the}
\eea
\end{conjecture}
Here it is essential to restrict to $2\times N$ partitionings as
for $3\times 3$ partitionings counterexamples can be found \cite{ka07}.

A special case of the conjectured inequality would imply validity of Hanner's inequality
for operators.
Hanner's inequality for $L_p$ functions $f,g$ is
\be
||f+g||_p^p + ||f-g||_p^p \ge (||f||_p+||g||_p)^p + \big|\, ||f||_p-||g||_p \big|^p,
\ee
for $1\le p\le 2$, and the reversed inequality for $2\le p$.
It is widely believed that these inequalities are also true for the Schatten trace ideals $\cS_p$:
for operators $A,B$ in $\cS_p$, and $1\le p\le 2$, that would mean
\be\label{hanner}
||A+B||_p^p + ||A-B||_p^p \ge (||A||_p+||B||_p)^p + \big|\, ||A||_p-||B||_p \big|^p,
\ee
while for $2\le p$, the inequality is reversed.
This generalisation has been proven in a number of instances:
\begin{enumerate}
\item For $1\le p$, when $A$ and $B$ are matrices such that $A+B$ and $A-B$ are positive semidefinite.
\item For $1\le p\le4/3$, $p=2$, and $4\le p\le\infty$, when $A$ and $B$ are general matrices.
\end{enumerate}
Proofs are due to Ball, Carlen and Lieb \cite{lieb} and Tomczak-Jaegermann \cite{TJ}.

That Conjecture \ref{conj_ka1} implies Hanner's inequality for matrices can be seen easily
by putting $N=2$, $A_1=B_2=A$ and $A_2=B_1=B$.
Unitarily conjugating
$T$ with the matrix $\frac{1}{\sqrt{2}}\twomat{\id}{\id}{\id}{-\id}$, and its norm compression with
$\frac{1}{\sqrt{2}}\twomat{1}{1}{1}{-1}$ directly yields (\ref{hanner}).

For $2\times 2$-partitioned block matrices $T$ ($N=2$) there are two other
special cases where Conjecture \ref{conj_ka1} is
known to hold, namely when $T$ is PSD (proven by King \cite{king}),
and when the blocks of $T$ are all diagonal matrices
(proven by King and Nathanson \cite{king_nath}).

Conjecture \ref{conj_ka1} has been proven for the following special cases \cite{ka07}:
the norm compression of $T$ has rank 1;
all blocks in $T$ have rank 1;
all blocks $A_k$ in the first row are proportional, and so are all blocks $B_k$ in the second row;
general $2\times N$-partitioned $T$, for $p\ge4$.
\section{A trace norm inequality for commutators}
The field of quantum information theory offers a rich source of problems in matrix theory.
In this section we consider one such problem that was raised by Bravyi \cite{bravyi} and
proven by him in a very special case.
\begin{conjecture}
Let $A,B$ be positive definite matrices such that $\trace(A+B)=1$; then
there is a constant $c$ (independent of the dimensions of $A$ and $B$) such that
\be
||\,\, [B,\log(A+B)]\,\,||_1 \le c (-\trace A \log \trace A -\trace B\log \trace B).
\ee
\end{conjecture}
Numerical work indicates that $c$ might actually be 1.

The more general case of $t:=\trace(A+B)\neq 1$ can be reduced to this inequality via
the transformation $A\mapsto A/t$ and $B\mapsto B/t$. After simplifying this yields
\bea
\lefteqn{||\,\, [B,\log(A+B)]\,\,||_1} \nonumber \\
& \le& c (\trace(A+B)\log \trace(A+B)-\trace A \log \trace A -\trace B\log \trace B),
\label{eq:SIM1}
\eea
which should hold for all $A,B\ge0$.
Note that the right-hand side is a function of $a:=\trace(A)$ and $b:=\trace(B)$ only.

This begs the question whether similar inequalities might be found when replacing the logarithm by
any function $f$; that is,
for a given function $f$, is there a constant $c$ and a function $g(a,b)$ such that
$$
||\,\,[B,f(A+B)]\,\,||_1 \le c\,g(a,b)?
$$
Since $A+B$ commutes with $f(A+B)$, the left-hand side can be replaced by
$||\,\,[A,f(A+B)]\,\,||_1$. Hence, $g(a,b)$ should be symmetric in its arguments.
Furthermore, for $a$ smaller than $b$, $g(a+b)$ should behave roughly like $a/(a+b)$.

The right-hand side of (\ref{eq:SIM1}) can be expressed in terms of the definite integral
$F(x)=\int_0^x \log(y)dy = x\log x -x$ as $F(a+b)-F(a)-F(b)$.
A preliminary numerical study indicated that this might be a useful way of looking at the
inequality (\ref{eq:SIM1}), leading to a more general conjecture:
\begin{conjecture}
Let $A$ and $B$ be positive semidefinite $d\times d$ matrices with $a=\trace A$ and $b=\trace B$.
For certain functions $f(x): \R\mapsto \R$ (of a class to be determined)
there exists a constant $c$, independent of $d$, such that
\be
||\,\, [B,f(A+B)]\,\,||_1 \le c (F(a+b)-F(a)-F(b)),
\label{eq:SIM2}
\ee
where $F(x)=\int_0^x f(y)dy$.
\end{conjecture}
Numerical experiments also revealed that $c=1$ is the best constant in (\ref{eq:SIM2})
for the functions $f(x)=\log(x)$, $f(x)=x^2$, $f(x)=x^p$ with $0<p\le 1$,
at least under the restriction $a+b=1$.

This is, in fact, very easy to prove for the functions $f(x)=x$ and $f(x)=x^2$.
Firstly, for $f(x)=x$, we get
$||\,\,[B,A+B]\,\,||_1=||\,\,[B,A]\,\,||_1\le ||A||_1 ||B||_1 =ab$
\cite{aw}, while $F(a+b)-F(a)-F(b)=((a+b)^2-a^2-b^2)/2=ab$. In this case we find that $c=1$.

Secondly, for $f(x)=x^2$, we get
\beas
||\,\,[B,(A+B)^2]\,\,||_1 &=& ||\,\,[B,A^2+AB+BA+B^2]\,\,||_1 \\
&=& ||\,\,([B,A^2]+[B^2,A])\,\,||_1 \\
&\le& ||\,\,[B,A^2]\,\,||_1+||\,\,[B^2,A]\,\,||_1\\
&\le& ba^2+b^2a = ((a+b)^3-a^3-b^3)/3.
\eeas
Here we exploited the fact that for $X,Y\ge0$, $||\,\,[X,Y]\,\,||_1 \le \trace(X)\trace(Y)$
and $\trace X^2\le (\trace X)^2$.
Again, we find that $c=1$.
\section{A majorisation relation}
Denote the singular values of a matrix $X$, arranged in non-increasing order,
by $\sigma_i(X)$.
Let $f$ be a concave function $f:\R_+\mapsto\R_+$ with $f(0)=0$,
then the following inequality holds:
\be
\sum_{i=1}^k (f(\sigma_i(X)) - f(\sigma_i(Y))) \le
\sum_{i=1}^k f(\sigma_i(X-Y)).\label{eq:uchi}
\ee
This is an immediate corollary of Theorem 4.4 in \cite{uchiyama}.

W.~Miao \cite{miao} conjectured that a stronger inequality holds,
whereby the sum in the left-hand side is replaced
by a sum of absolute values:
\begin{conjecture}[Miao]
Let $X,Y\in\M_{n,m}(\C)$.
Let $f$ be a concave function $f:\R_+\mapsto\R_+$ with $f(0)=0$.
Then for any $k\le n,m$,
\be
\sum_{i=1}^k |f(\sigma_i(X)) - f(\sigma_i(Y))| \le
\sum_{i=1}^k f(\sigma_i(X-Y)).\label{eq:0}
\ee
\end{conjecture}
The particular case of this inequality for $f(x)=x^q$, $0<q<1$,
would have an application in the study of compressed sensing methods \cite{oymak}.

The conjecture is relatively easy to prove when $X$ and $Y$ are diagonal matrices; we show this below.

We denote by $\mathbf{x}^\u$ a rearrangement of the vector $\mathbf{x}$ in non-decreasing order, and by
$\mathbf{x}^\d$ a rearrangement in non-increasing order.

Let $f$ be a concave function $f:\R_+\mapsto\R_+$ with $f(0)=0$.
It follows that $f$ is non-decreasing and subadditive, i.e.\ $f(x)\le f(x+y)\le f(x)+f(y)$ for $x,y\ge0$.

\begin{proposition}
Let $\mathbf{u}=(u_1,u_2,\ldots,u_n)$ and $\mathbf{v}=(v_1,v_2,\ldots,v_n)$, with $u_i,v_i\in\C$.
Let $\mathbf{x}=|\mathbf{u}|$, $\mathbf{y}=|\mathbf{v}|$
and $\mathbf{z}=|\mathbf{u}-\mathbf{v}|$, i.e.\ $z_i=|u_i-v_i|$.
For any non-negative, monotonically increasing, subadditive function $f$
\be
\sum_{i=1}^k |f(x_i^\d)-f(y_i^\d)| \le
\sum_{i=1}^k f(z_i^\d),\quad k=1,\ldots,n. \label{eq:1}
\ee
\end{proposition}
\textit{Proof.}
Let $\mathbf{a}=f(\mathbf{x})$ and $\mathbf{b}=f(\mathbf{y})$.
Because $f$ is non-decreasing, $\mathbf{a}^\d=f(\mathbf{x}^\d)$,
$\mathbf{b}^\d=f(\mathbf{y}^\d)$ and $f(\mathbf{z})^\d=f(\mathbf{z}^\d)$.
Let $\mathbf{c}=|\mathbf{a}^\d - \mathbf{b}^\d|$ and $\mathbf{d}=|\mathbf{a}-\mathbf{b}|$.
The LHS of (\ref{eq:1}) is thus given by $\sum_{i=1}^k c_i$.
Obviously, we have
\be
\mbox{LHS} = \sum_{i=1}^k c_i \le \sum_{i=1}^k c_i^\d.
\ee
Proposition 6.A.2.a in \cite{MOA} states $|\mathbf{a}^\d - \mathbf{b}^\d| \prec_w |\mathbf{a}-\mathbf{b}|$.
Thus, $\mathbf{c}\prec_w\mathbf{d}$, i.e.
\be
\sum_{i=1}^k c_i^\d \le \sum_{i=1}^k d^\d_i.
\ee

Because $f$ is subadditive, for $x\ge y\ge0$ we have $f(x)-f(y)\le f(x-y)$.
For $y\ge x\ge0$ we obviously have $f(y)-f(x)\le f(y-x)$.
More generally, then,
for $x,y\ge0$ we have $|f(x)-f(y)|\le f(|x-y|)$.

Thus, $d_i= |a_i-b_i| =  |f(x_i)-f(y_i)| \le f(|x_i-y_i|)$.
Furthermore, since $||u|-|v||\le |u-v|$ for $u,v\in\C$, we have
$f(|x_i-y_i|) = f(||u_i|-|v_i||)\le f(|u_i-v_i|)= f(z_i)$. In other words
$\mathbf{d} \le f(\mathbf{z})$, which implies
$\mathbf{d} \prec_w f(\mathbf{z})$.
We obtain
\be
\sum_{i=1}^k d^\d_i \le \sum_{i=1}^k f(\mathbf{z})^\d_i = \sum_{i=1}^k f(z_i^\d) = \mbox{RHS}.
\ee
Combining the three inequalities yields (\ref{eq:1}).
In fact this proves the slightly stronger statement
\be
|f(\mathbf{|u|}^\d)-f(\mathbf{|v|}^\d)| \prec_w f(|\mathbf{u}-\mathbf{v}|). \label{eq:1b}
\ee
\qed

The conjecture can now be proven easily in the simple case that $X$ and $Y$ are complex diagonal matrices.
Indeed, let $X=\diag(\mathbf{u})$ and $Y=\diag(\mathbf{v})$,
then $\sigma(X)=\mathbf{|u|}^\d$, $\sigma(Y)=\mathbf{|v|}^\d$
and $\sigma(X-Y) = \mathbf{z}^\d$, with $\mathbf{z}=|\mathbf{u}-\mathbf{v}|$.
The proposition immediately yields the required statement.
\section{Inequalities for several operators}
\subsection{Clarkson inequalities for several operators}
The classical Clarkson inequalities for two operators $A$ and $B$ on a
separable Hilbert space assert that%
\begin{equation}
2\left( \left\Vert A\right\Vert _{p}^{p}+\left\Vert B\right\Vert
_{p}^{p}\right) \leq \left\Vert A+B\right\Vert _{p}^{p}+\left\Vert
A-B\right\Vert _{p}^{p}\leq 2^{p-1}\left( \left\Vert A\right\Vert
_{p}^{p}+\left\Vert B\right\Vert _{p}^{p}\right)  \label{Eq1}
\end{equation}%
for $2\leq p<\infty $,
\begin{equation}
2^{p-1}\left( \left\Vert A\right\Vert _{p}^{p}+\left\Vert B\right\Vert
_{p}^{p}\right) \leq \left\Vert A+B\right\Vert _{p}^{p}+\left\Vert
A-B\right\Vert _{p}^{p}\leq 2\left( \left\Vert A\right\Vert
_{p}^{p}+\left\Vert B\right\Vert _{p}^{p}\right)  \label{Eq2}
\end{equation}%
for $0<p\leq 2$,%
\begin{equation}
2\left( \left\Vert A\right\Vert _{p}^{p}+\left\Vert B\right\Vert
_{p}^{p}\right) ^{q/p}\leq \left\Vert A+B\right\Vert _{p}^{q}+\left\Vert
A-B\right\Vert _{p}^{q}  \label{Eq3}
\end{equation}%
for $2\leq p<\infty $; $\frac{1}{p}+\frac{1}{q}=1$, and
\begin{equation}
\left\Vert A+B\right\Vert _{p}^{q}+\left\Vert A-B\right\Vert _{p}^{q}\leq
2\left( \left\Vert A\right\Vert _{p}^{p}+\left\Vert B\right\Vert
_{p}^{p}\right) ^{q/p}  \label{Eq4}
\end{equation}%
for $1<p\leq 2$; $\frac{1}{p}+\frac{1}{q}=1.$

These inequalities have been generalized to $n$-tuples of operators by
Bhatia and Kittaneh \cite{BK} as follows:

If $A_{0},A_{1},\cdot \cdot \cdot ,A_{n-1}$ are operators on a separable
Hilbert space, and $\omega _{0},\omega _{1},\cdot \cdot \cdot ,\omega _{n-1}$
are the $n$-th roots of unity with $\omega _{j}=e^{2\pi ij/n},\,\,0\leq j\leq
n-1$, then for every unitarily invariant norm,
\begin{equation}
n\sum_{j=0}^{n-1}\left\Vert A_{j}\right\Vert _{p}^{p}\leq
\sum_{k=0}^{n-1}\left\Vert \sum_{j=0}^{n-1}\omega _{j}^{k}A_{j}\right\Vert
_{p}^{p}\leq n^{p-1}\sum_{j=0}^{n-1}\left\Vert A_{j}\right\Vert _{p}^{p}
\label{Eq5}
\end{equation}%
for $2\leq p<\infty $,
\begin{equation}
n^{p-1}\sum_{j=0}^{n-1}\left\Vert A_{j}\right\Vert _{p}^{p}\leq
\sum_{k=0}^{n-1}\left\Vert \sum_{j=0}^{n-1}\omega _{j}^{k}A_{j}\right\Vert
_{p}^{p}\leq n\sum_{j=0}^{n-1}\left\Vert A_{j}\right\Vert _{p}^{p}
\label{Eq6}
\end{equation}%
for $0<p\leq 2$,%
\begin{equation}
n\left( \sum_{j=0}^{n-1}\left\Vert A_{j}\right\Vert _{p}^{p}\right)
^{q/p}\leq \sum_{k=0}^{n-1}\left\Vert \sum_{j=0}^{n-1}\omega
_{j}^{k}A_{j}\right\Vert _{p}^{q}  \label{Eq7}
\end{equation}%
for $2\leq p<\infty $; $\frac{1}{p}+\frac{1}{q}=1$, and
\begin{equation}
\sum_{k=0}^{n-1}\left\Vert \sum_{j=0}^{n-1}\omega _{j}^{k}A_{j}\right\Vert
_{p}^{q}\leq n\left( \sum_{j=0}^{n-1}\left\Vert A_{j}\right\Vert
_{p}^{p}\right) ^{q/p}  \label{Eq8}
\end{equation}%
for $1<p\leq 2$; $\frac{1}{p}+\frac{1}{q}=1$.

Other natural generalizations of the inequalities (\ref{Eq1}) and (\ref{Eq2})
to $n$-tuples of operators have been recently given by Hirzallah and
Kittaneh \cite{HK}:

If $A_{0},A_{1},\cdot \cdot \cdot ,A_{n-1}$ are operators on a separable
Hilbert space, then
\begin{equation}
\left\Vert \sum_{j=0}^{n-1}A_{j}\right\Vert _{p}^{p}+\sum_{0\leq j<k\leq
n-1}\left\Vert A_{j}-A_{k}\right\Vert _{p}^{p}\leq
n^{p-1}\sum_{j=0}^{n-1}\left\Vert A_{j}\right\Vert _{p}^{p}  \label{Eq9}
\end{equation}
for $2\leq p<\infty $ and
\begin{equation}
n^{p-1}\sum_{j=0}^{n-1}\left\Vert A_{j}\right\Vert _{p}^{p}\leq \left\Vert
\sum_{j=0}^{n-1}A_{j}\right\Vert _{p}^{p}+\sum_{0\leq j<k\leq n-1}\left\Vert
A_{j}-A_{k}\right\Vert _{p}^{p}  \label{Eq10}
\end{equation}
for $0<p\leq 2$.

In view of the inequalities (\ref{Eq9}) and (\ref{Eq10}), it is reasonable
to make the following conjecture concerning new natural generalizations
(along the lines of the inequalities (\ref{Eq9}) and (\ref{Eq10})) of the
inequalities (\ref{Eq3}) and (\ref{Eq4}) to $n$-tuples of operators.

\begin{conjecture}
Let $A_{0},A_{1},\cdot \cdot \cdot ,A_{n-1}$ be operators on a separable
Hilbert space. Then%
\begin{equation}
n\left( \sum_{j=0}^{n-1}\left\Vert A_{j}\right\Vert _{p}^{p}\right)
^{q/p}\leq \left\Vert \sum_{j=0}^{n-1}A_{j}\right\Vert _{p}^{q}+\sum_{0\leq
j<k\leq n-1}\left\Vert A_{j}-A_{k}\right\Vert _{p}^{q}  \label{Eq11}
\end{equation}%
for $2\leq p<\infty $; $\frac{1}{p}+\frac{1}{q}=1$ and
\begin{equation}
\left\Vert \sum_{j=0}^{n-1}A_{j}\right\Vert _{p}^{q}+\sum_{0\leq j<k\leq
n-1}\left\Vert A_{j}-A_{k}\right\Vert _{p}^{q}\leq n\left(
\sum_{j=0}^{n-1}\left\Vert A_{j}\right\Vert _{p}^{p}\right) ^{q/p}
\label{Eq12}
\end{equation}%
for $0<p\leq 2$; $\frac{1}{p}+\frac{1}{q}=1$.
\end{conjecture}
\subsection{Hlawka-type inequalities for operators}
The extension of the triangle inequality to several operators is immediate and of course well-known.
Let $X$, $Y$ and $Z$ be arbitrary matrices in $\M_{N,M}(\C)$,
 $||.||$ be a norm on $\M_{N,M}(\C)$.
Then the triangle inequality reads $||X||+||Y||+||Z|| - ||X+Y+Z||\ge 0$.

Highly non-trivial questions are obtained when we ask for  complementary inequalities,
as we did in Section \ref{sec:compl} for the McCarthy inequality.
Here we consider an upper bound on $||X||+||Y||+||Z|| - ||X+Y+Z||$  in
terms of the pairwise quantities $||X||+||Y||-||X+Y||$, $||X||+||Z||-||X+Z||$
and $||Y||+||Z||-||Y+Z||$.

\begin{problem}
Find the best (dimension-independent) constant $\cC$ -- if it exists -- such that
\bea
\lefteqn{||X||+||Y||+||Z|| - ||X+Y+Z||} \nonumber\\
&\le& \cC \{(||X||+||Y||-||X+Y||)+\nonumber\\
&&(||X||+||Z||-||X+Z||)+\nonumber\\
&&(||Y||+||Z||-||Y+Z||)
\}\label{eq:hlawka}
\eea
In particular, find the best constant $\cC_p$ when the norm in question is the Schatten $p$-norm.
\end{problem}
There are several reasons to study this question.
First of all, it has a nice geometrical interpretation.
Consider the norm ball $\cB$ of $||\cdot||$.
For the sake of illustration, consider only $X,Y,Z$ of norm $1$,
i.e.\ lying on the surface of the norm ball.
The matrices $X$, $Y$ and $Z$ support a 2D plane:
$$
\cP = \{xX+yY+zZ: x,y,z\in\R, x+y+z=1\}
$$
Inequality (\ref{eq:hlawka}) tells us something about the shape of the intersection $\cB\cap\cP$.
Assume that $\cP$ is not tangent to $\cB$, i.e.\ $X$, $Y$ and $Z$ do not lie on a facet.
Define $f(x,y,z) = x||X||+y||Y||+z||Z|| - ||xX+yY+zZ||$.
Thus $f(1,0,0)=f(0,1,0)=f(0,0,1)=0$.
Then (\ref{eq:hlawka}) becomes
$$
f(1/3,1/3,1/3) \le \frac{2\cC}{3}(f(1/2,1/2,0) + f(1/2,0,1/2)+f(0,1/2,1/2)).
$$
This says that the intersection \textit{can never be a triangle}.
Otherwise, we would have
$$
f(1/2,1/2,0)=f(1/2,0,1/2)=f(0,1/2,1/2)=0,
$$
which cannot be
unless $f(1/3,1/3,1/3)=0$, contradicting the assumption that $\cP$ is not tangent to $\cB$.

For the infinity norm $||\cdot||_\infty$, it turns out that $\cC_\infty$ is infinite.
If we take $X=e^{11}+e^{22}$, $Y=e^{11}+e^{33}$, $Z=e^{22}+e^{33}$, then we get
\beas
||X||+||Y||+||Z|| - ||X+Y+Z|| &=& 1, \\
\mbox{whereas }||X||+||Y||-||X+Y||&=&0, \mbox{ etc } \ldots
\eeas
This is clear geometrically: consider vectors in $\R^3$, so that the norm ball is a cube.
There are of course infinitely many planes whose intersection with a cube is a triangle.

For the trace norm, numerical simulations show that $\cC_1$ is at least $40$, indicating
that it might be arbitrarily large as well. This is relevant for a conjectured inequality
mentioned in the next section.

For the $2$-norm, the sharpest value of $\cC_2$ is known to be $1$, as (\ref{eq:hlawka})
then reduces to Hlawka's inequality \cite{mitr},
which was originally phrased for vectors in $\R^n$ as:
\be
||\vec{x}+\vec{y}+\vec{z}|| + ||\vec{x}||+||\vec{y}||+||\vec{z}||
\ge ||\vec{x}+\vec{y}|| + ||\vec{x}+\vec{z}|| + ||\vec{y}+\vec{z}||.
\ee
Here, the norm is any vector norm that satisfies the parallellogram identity.
Hence, it also holds for matrices in the Hilbert-Schmidt norm $||\cdot||_2$.

There is a generalisation of Hlawka's inequality to $K$ vectors,
due to Adamovic \cite{mitr}:
\be
||\sum_{j=1}^K \vec{x_j}|| +(K-2)\sum_{j=1}^K ||\vec{x_j}|| \ge \sum_{1\le j<k\le K} ||\vec{x_j}+\vec{x_k}||.
\ee
Thus, once the problem of finding the best constant in (\ref{eq:hlawka}) has been tackled
one might consider generalisations of (\ref{eq:hlawka}) to any number of matrices.

\subsection{An inequality with relevance for quantum hypothesis testing}
Related to both conjectured multivariate inequalities presented above is the
following:
\begin{conjecture}
Let $A_j$ be $n$ positive operators, then
\be
\trace\sqrt{\sum_j A_j^2} - \sum_j\trace A_j
\ge c \sum_{j<k} (||A_j-A_k||_1 - \trace( A_j + A_k)).\label{eq:47}
\ee
\end{conjecture}
where $c$ is a constant.

This question is of interest in quantum information science and would imply an
important new result in quantum hypothesis testing. The actual value of $c$ is
of no importance for this, even though we believe it to be $1$.

A weaker, but more symmetric looking block matrix inequality that would imply the previous
inequality is:
\be
||(A_1\;\cdots\;A_n)||_1 - \sum_j\trace A_j
\ge c' \sum_{j<k} (||(A_j\;\;A_k)||_1 - \trace (A_j+ A_k)).
\ee
Numerical calculations point at a value of $c'$ of roughly $2.7$, for $n=3$.
Note that this inequality would follow from the Hlawka inequality for the trace norm
if $\cC_1$ were finite (contrary to available numerical evidence).

To obtain the previous inequality from this, we use the inequality
$$
\trace(A+B)-||(A\;\; B)||_1 = \trace(A+B)-\trace\sqrt{A^2+B^2}
\le(\trace(A+B)-||A-B||_1)/2,
$$
which is easy to prove from concavity of the function $\trace\sqrt{X}$.

A number of further explorations are possible. One could envisage replacing the term
$\trace\sqrt{\sum_j A_j^2}$ by $\trace(\sum_j A_j^{2p})^{1/2p}$ and
the traces and trace norm in the other terms by Schatten $p$-norms. Furthermore, can
the inequality (\ref{eq:47}) be modified so that it holds for general (non-positive, non-normal) operators?

\section{Geometry of polynomials}
Broadly speaking, the geometry of polynomials is the study of the location in the complex plane $\C$
of the roots of a polynomial $p(z)$ and how this relates to other properties of the polynomial, such as
the location of its coefficients in $\C$, or the location of the roots of its derivative $p'(z)$.

Here we will concern ourselves with the latter.
We shall denote the roots of the degree $n$ polynomial $p(z)$ by $\lambda_i$, $i=1,\ldots,n$,
and the roots of $p'(z)$ by $\mu_j$, $j=1,\ldots,n-1$.
By inspection of the coefficient of $p(z)$ of degree $n-1$ and the coefficient of $p'(z)$ of degree $n-2$
it is easy to see that both sets of roots have the same arithmetic mean, a quantity which we'll denote by $G$:
\be
G=\frac{1}{n}\sum_{i=1}^n \lambda_i = \frac{1}{n-1}\sum_{j=1}^{n-1}\mu_i.
\ee

One of the first, and easiest, results relating these two sets of roots is the Gauss-Lucas theorem,
according to which all $\mu_j$ are in the convex hull of the set $\{\lambda_i\}_{i=1}^n$.
A great many relations have been discovered that improve on this.
For example, according to Schoenberg's conjecture, which has been proven recently \cite{malamud,pereira},
\be
\sum_{j=1}^{n-1}|\mu_j|^2 \le |G|^2 + \frac{n-2}{n}\,\sum_{i=1}^n |\lambda_i|^2.
\ee
Many other relations have been conjectured (see e.g. \cite{khavinson,schmeisser}):
\begin{conjecture}[Smale]
Let $p$ be a polynomial of degree $n$ such that $p(0)=0$ and $p'(0)\neq 0$.
At least one of the roots of $p'$ satisfies $p(\mu_j)/\mu_j \le p'(0)$.
\end{conjecture}

The following conjecture has been open since 1959:
\begin{conjecture}[Sendov]
Let $p(z) = (z- \lambda_1)\cdots
(z-\lambda_n)$ be such that $|\lambda_i| \le 1$ for all $i$, and
$p^{\prime}(z) = n(z-\mu_1)\cdots (z-\mu_{n-1})$. Then  for any $\lambda_i$
there is $\mu_j$ such that $|\lambda_i - \mu_j| \le 1$.
\end{conjecture}
This conjecture has also been known as ``Ilyeff's conjecture'' (as in, e.g. \cite{rubinstein}); see \cite{schmeisser}
for the history of the problem.

Matrix analysis is highly relevant to the study of the geometry of polynomials.
Polynomials and matrices are, indeed, intimately related: the eigenvalues of a matrix are the roots
of its characteristic polynomial, and, conversely, for any polynomial there is a matrix,
called the companion matrix, having that polynomial as characteristic polynomial.
For example, Pereira's proof of Schoenberg's conjecture proceeds by reducing it to a matrix inequality involving the
sums of the squares of the absolute values of the eigenvalues of a normal matrix.

Cheung and Ng \cite{cheung}
introduced the concept of $D$-companion matrix, which they used to prove a conjecture by de Bruijn and Sharma.
Let $I$ and $J$ be the identity matrix of order $n-1$ and the $(n-1)\times (n-1)$ matrix with all entries equal
to 1, respectively, then the $D$-companion matrix of the degree $n$ polynomial $p$ with roots $\lambda_i$
is the matrix $\diag(\lambda_1,\ldots,\lambda_{n-1})(I-J/n)+\lambda_n J/n$. Cheung and Ng showed
that its eigenvalues are exactly the roots $\mu_j$ of $p'$.

This leads to the following reformulation of Sendov's conjecture:
\begin{conjecture}
Let $D = \diag(\lambda_1 - \lambda_n,
\dots, \lambda_{n-1}-\lambda_n)$, where $\lambda_1, \ldots, \lambda_n$ belong to the closed unit disk.
Then $D(I - J/n)$ has an eigenvalue in the closed unit disk.
\end{conjecture}
For other applications of the $D$-companion matrix to polynomial geometry, see e.g. \cite{adm}.

Another reformulation of Sendov's conjecture exploits certain properties of circulant matrices.
The circulant matrix $C$ with vector of coefficients $\mathbf{c}=(c_0,\ldots,c_{n-1})$, denoted
$C=\cirk(c_0,\ldots,c_{n-1})$, is the matrix whose first row is given by the vector $\mathbf{c}$
and whose subsequent rows are obtained form the first one by successive cyclic permutations;
i.e.
$$
C=\cirk(c_0,\ldots,c_{n-1}) = \left(
\begin{array}{ccccc}
c_0 & c_1 & \cdots & c_{n-2} & c_{n-1} \\
c_{n-1} & c_0 & c_1 & & c_{n-2} \\
\vdots & c_{n-1} & c_0& \ddots & \vdots \\
c_2 & & \ddots& \ddots & c_1 \\
c_1 & c_2 & \cdots & c_{n-1} & c_0
\end{array}
\right).
$$
Because of its inherent symmetry, the principal submatrices of $C$ obtained by deleting a single row
and its corresponding column are identical, up to conjugation with a permutation matrix.
Let $C'$ be such a principal submatrix.

The characteristic polynomial of $C$ is given by
$$
p(z)=\det(C-z\id)=\det\cirk(c_0-z,c_1,\ldots,c_{n-1}).
$$
Using the abovementioned property of the principal submatrices, it is easy to see that the derivative
of $p(z)$ is given by $p'(z)=-n\det(C'-z\id)$. In other words, the roots of $p$ are the
eigenvalues of $C$ and the roots of $p'$ are the eigenvalues of $C'$.

This yields a third formulation of Sendov's conjecture:
\begin{conjecture}
Suppose $C \in M_n$ is a
circulant matrix with $\|C\| \le 1$ and $\lambda$ is an eigenvalue of $C$.
Let $C^{\prime}$ be the principal submatrix obtained from $C$ by removing the
first row and the first column. Then $C^{\prime}-\lambda I$ has an eigenvalue in
the closed unit disk.
\end{conjecture}

\section{Miscellaneous Problems}
To end this chapter, we pose two ``inverse eigenvalue problems'':
\begin{problem} Determine the necessary and sufficient
conditions on $a_1 \ge \cdots \ge a_n$, $b_1 \ge \cdots \ge b_n$ and $s_1
\ge \cdots \ge s_n$ for the existence of Hermitian matrices $A$ and $B$ so
that $A$ has eigenvalues $a_1, \dots, a_n$, $B$ has eigenvalues $b_1, \dots,
b_n$, and $A+iB$ has singular values $s_1, \dots, s_n$.
\end{problem}

\begin{problem}
Determine the necessary and sufficient
conditions on the numbers $a_1, \ldots, a_k$ and $b_1, \ldots, b_n$ for the existence of
normal matrices $A \in M_k$ and $B \in M_n$ so that the matrices have these
numbers as eigenvalues and $A$ is a compression of $B$.
\end{problem}

Finally, we mention one of the problems put forward by Michel Crouzeix in \cite{crouzeix}
concerning a functional calculus based on the numerical range.
Let $W(A)$ be the numerical range of $A$ and $\|A\|$ the spectral
norm of $A$.
\begin{problem}[Crouzeix]
Let $p(z)$ be a polynomial. Find
the smallest $k$ such that for any matrix or operator $A$
\[
\|p(A)\| \le k \sup\{ |p(\mu)|: \mu \in W(A)\}.
\]
\end{problem}


\begin{thebibliography}{99}
\bibitem{adm} M.~Adm and F.~Kittaneh,
\textit{Bounds and majorization relations for the critical points of polynomials},
Linear Algebra Appl.\ \textbf{436} (2012), 2494--2503.

\bibitem{ka05} K.M.R. Audenaert, \textit{A norm compression inequality for
block partitioned positive definite matrices}, Linear Algebra
Appl.\ \textbf{413} (2006), 155--176.

\bibitem{ka06} K.M.R. Audenaert, \textit{A singular value for Heinz means},
Linear Algebra Appl.\ \textbf{422} (2007), 279--283.

\bibitem{ka07} K.M.R. Audenaert, \textit{On a norm compression inequality for $2\times N$ partitioned
block matrices}, Linear Algebra Appl.\ \textbf{428} (2008), 781--795.

\bibitem{ka2} K.M.R. Audenaert, \textit{On the Araki-Lieb-Thirring inequality},
Int.\ J.\ Inf.\ Syst.\ Sci.\ \textbf{4} (2008), 78--83.

\bibitem{ka_ilas} K.M.R. Audenaert, \textit{In-betweenness, a geometrical monotonicity
property for operator means}, Linear Algebra Appl., in press, 2012.

\bibitem{katrace} K.M.R.~Audenaert,
\textit{Trace inequalities for completely monotone functions and Bernstein functions},
Linear Algebra Appl., in press, 2012. See also eprint arXiv:1109.3057.

\bibitem{AB} J.S.~Aujla and J.-C.~Bourin, \textit{Eigenvalue inequalities for convex and log-convex functions},
Linear Algebra Appl.\ \textbf{424} (2007), 25--35.

\bibitem{lieb} K.~Ball, E.~Carlen and E.H.~Lieb, \textit{Sharp uniform convexity and smoothness
inequalities for trace norms}, Invent.\ Math.\ \textbf{115} (1994), 463--482.

\bibitem{bhagwat} K.V.~Bhagwat and A.~Subramanian, \textit{Inequalities between means of positive operators},
Math.\ Proc.\ Camb.\ Phil.\ Soc.\ \textbf{83} (1978), 393--401.

\bibitem{bhatia} R.~Bhatia, \textit{Matrix Analysis}, Springer, Heidelberg, 1997.

\bibitem{bk90} R.~Bhatia and F.~Kittaneh, \textit{Norm inequalities for partitioned
operators and an application}, Math.\ Ann.\ \textbf{287} (1990), 719--726.

\bibitem{BK} R.~Bhatia and F.~Kittaneh, \textit{Clarkson inequalities with several
operators}, Bull.\ London Math.\ Soc.\ \textbf{36} (2004), 820--832.

\bibitem{B} J.-C.~Bourin, \textit{Matrix subadditivity inequalities and
block-matrices}, Internat.\ J.\ Math.\ \textbf{20} (2009), 679--691.

\bibitem{B10} J.-C.~Bourin, \textit{A matrix subadditivity inequality for symmetric norms},
Proc.\ AMS \textbf{138} (2010), 495--504.

\bibitem{bhl} J.-C.~Bourin, T.~Harada and E.-Y.~Lee, \textit{Subadditivity inequalities for compact operators},
To appear in Canadian Bull.\ Math.\ (2012).

\bibitem{BU} J.-C.~Bourin and M.~Uchiyama, \textit{A matrix subadditivity inequality
for $f(A+B)$ and $f(A)+f(B)$}, Linear Algebra Appl.\ \textbf{423} (2007), 512--518.

\bibitem{bravyi} S.~Bravyi, \textit{Upper bounds on entangling rates of bipartite Hamiltonians},
Phys.\ Rev.\ A \textbf{76} (2007), 052319.

\bibitem{cheung} W.S.~Cheung and T.W.~Ng, \textit{A companion matrix approach to the study of zeros and
critical points of a polynomial}, J.\ Math.\ Anal.\ Appl.\ \textbf{319} (2006), 690--707.

\bibitem{crouzeix} M.~Crouzeix, \textit{Open problems on numerical range and functional calculus},
unpublished, 2006.
Downloaded 19/12/2011 from the author's personal web page {\tt http://perso.univ-rennes1.fr/michel.crouzeix}.

\bibitem{hanner} O.~Hanner, \textit{On the uniform convexity of $L^p$ and $l^p$},
Ark.\ Math.\ \textbf{3} (1958), 239--244.

\bibitem{HK} O. Hirzallah and F. Kittaneh, \textit{Non-commutative Clarkson
inequalities for $n$-tuples of operators}, Integral Equations Operator Theory
\textbf{60} (2008), 369--379.

\bibitem{HJI} R.A.~Horn and C.~R.~Johnson, \textit{Matrix Analysis}, Cambridge University Press, Cambridge, 1985.

\bibitem{khavinson} D.~Khavinson, R.~Pereira, M.~Putinar, E.B.~Saff and S.~Shimoron,
\textit{Borcea's variance conjectures on the critical points of polynomials},
in: Notions of Positivity and the Geometry of Polynomials, P.
Br\"and\'en et al. (eds.), Trends in Mathematics 2011, Birkh\"auser, Basel, 2011, 283--309.


\bibitem{king} C.~King, \textit{Inequalities for trace norms of $2\times 2$ block matrices},
Commun.\ Math.\ Phys.\ \textbf{242} (2003), 531--545.

\bibitem{king_nath} C.~King and M.~Nathanson, \textit{New trace norm inequalities for $2\times 2$
blocks of diagonal matrices},
Linear Algebra Appl.\ \textbf{389} (2004), 77--93.

\bibitem{K} T. Kosem, \textit{Inequalities between $\left\Vert f(A+B)\right\Vert $
and $\left\Vert f(A)+f(B)\right\Vert $}, Linear Algebra Appl.\ \textbf{418} (2006),
153--160.

\bibitem{kuboando} F.~Kubo and T.~Ando, \textit{Means of positive linear operators},
Math.\ Ann.\ \textbf{246} (1980), 205--224.


\bibitem{lee11} E.-Y.~Lee, \textit{How to compare the absolute values of operator sums and the sums of absolute values?},
In press, Operators and Matrices, 2012.

\bibitem{malamud} S.M.~Malamud, \textit{Inverse spectral problem for normal matrices and the Gauss-Lucas theorem},
Trans.\ Amer.\ Math.\ Soc.\ \textbf{357} (2005), 4043--4064.

\bibitem{miao} Weimin Miao, personal communication, 2011.

\bibitem{mitr} D.S.~Mitrinovic, J.E.~Pecaric and A.M.~Fink,
\textit{Classical and New Inequalities in Analysis}, Kluwer, Dordrecht, 1993.

\bibitem{MOA} A.W.~Marshall, I.~Olkin and B.C.~Arnold,
\textit{Inequalities: theory of majorization and its applications}, Second edition, Springer, Heidelberg, 2011.

\bibitem{oymak} S.~Oymak, K.~Mohan, M.~Fazel and B.~Hassibi,
\textit{A simplified approach to recovery conditions for low rank matrices},
Proc.\ Intl.\ Symp.\ Info.\ Theory (ISIT), August 2011; see also eprint arXiv:1103.1178 (2011).

\bibitem{pereira} R.~Pereira, \textit{Differentiators and the geometry of polynomials},
J.\ Math.\ Anal. Appl.\ \textbf{285} (2003), 336--348.

\bibitem{rubinstein} Z.~Rubinstein, \textit{On a problem of Ilyeff}, Pacific J.\ Math.\ \textbf{36} (1968), 159--161.

\bibitem{schmeisser} G.~Schmeisser, \textit{The conjectures of Sendov and Smale},
in: Approximation theory: A volume dedicated to Blagovest Sendov,
B.\ Bojanov (ed.), DARBA, Sofia, 2002, 353--369.

\bibitem{TJ} N.~Tomczak-Jaegermann, \textit{The moduli of smoothness and convexity and
Rademacher averages of trace classes $S_p$},
Studia Math.\ \textbf{50} (1974), 163--182.

\bibitem{uchiyama} M.~Uchiyama, \textit{Subadditivity of eigenvalue sums},
Proc.\ AMS \textbf{134} (2006), 1405--1412.

\bibitem{aw} D.~Wenzel and K.M.R.~Audenaert, \textit{Impressions of convexity --
An illustration for commutator bounds}, Linear Algebra Appl.\ \textbf{433} (2010), 1726--1759.

\end{thebibliography}
\end{document}